\documentstyle[12pt]{article}
\begin{document}
\renewcommand{\thefootnote}{\fnsymbol{footnote}}
\newpage
\pagestyle{empty}
\setcounter{page}{0}
\renewcommand{\thesection}{\arabic{section}}
\renewcommand{\theequation}{\thesection.\arabic{equation}}
\newcommand{\sect}[1]{\setcounter{equation}{0}\section{#1}}

%
%
\def\CC{{\Bbb C}}
\def\NN{{\Bbb N}}
\def\QQ{{\Bbb Q}}
\def\RR{{\Bbb R}}
\def\ZZ{{\Bbb Z}}
\def\cA{{\cal A}}          \def\cB{{\cal B}}          \def\cC{{\cal C}}
\def\cD{{\cal D}}          \def\cE{{\cal E}}          \def\cF{{\cal F}}
\def\cG{{\cal G}}          \def\cH{{\cal H}}          \def\cI{{\cal I}}
\def\cJ{{\cal J}}          \def\cK{{\cal K}}          \def\cL{{\cal L}}
\def\cM{{\cal M}}          \def\cN{{\cal N}}          \def\cO{{\cal O}}
\def\cP{{\cal P}}          \def\cQ{{\cal Q}}          \def\cR{{\cal R}}
\def\cS{{\cal S}}          \def\cT{{\cal T}}          \def\cU{{\cal U}}
\def\cV{{\cal V}}          \def\cW{{\cal W}}          \def\cX{{\cal X}}
\def\cY{{\cal Y}}          \def\cZ{{\cal Z}}
\def\qed{\hfill \rule{5pt}{5pt}}
\def\id{\mbox{id}}
\def\ggo{{\frak g}_{\bar 0}}
\def\uqggo{\cU_q({\frak g}_{\bar 0})}
\def\uqggp{\cU_q({\frak g}_+)}
\def\half{\frac{1}{2}}
\def\btf{\bigtriangleup}
\newtheorem{lemma}{Lemma}
\newtheorem{prop}{Proposition}
\newtheorem{theo}{Theorem}
\newtheorem{Defi}{Definition}

\vfill
\vfill
\begin{center}

{\LARGE {\bf {\sf
RTT relations, a modified braid equation and noncommutative planes
}}} \\[0.8cm]

{\large  A.Chakrabarti\footnote{chakra@cpht.polytechnique.fr}},

\begin{center}
{\em
Centre de Physique Th\'eorique\footnote{Laboratoire Propre
du CNRS UPR A.0014}, Ecole Polytechnique, 91128 Palaiseau Cedex, France.\\
}
\end{center}

\end{center}

\smallskip

\smallskip

\smallskip

\smallskip

\smallskip

\smallskip

\begin{abstract}
With the known group relations for the elements $(a,b,c,d)$ of a quantum
matrix $T$ as input a general solution of the $RTT$ relations is sought
without imposing the Yang - Baxter constraint for $R$ or the braid
equation for $\hat{R} = PR$. For three biparametric deformatios,
$GL_{(p,q)}(2), GL_{(g,h)}(2)$ and $GL_{(q,h)}(1/1)$, the standard,the
nonstandard and the hybrid one respectively, $R$ or $\hat{R}$ is found to
depend , apart from the two parameters defining the deformation in
question, on an extra free parameter $K$,such that

$\hat{R}_{(12)}\hat{R}_{(23)}\hat{R}_{(12)} -
\hat{R}_{(23)}\hat{R}_{(12)}\hat{R}_{(23)} = \biggl( { K\over {
K_{1}}} - 1
\biggr) \biggl( {K\over { K_{2}}} - 1 \biggr) \biggl(\hat{R}_{(23)} -
\hat{R}_{(12)}\biggr)$

with $ (K_{1},K_{2}) = (1,{p\over {q}}), (1,1)$ and $(1,{1\over {q}})$
respectively. Only for $K =K_{1}$ or $K= K_{2}$ one has the braid equation.
Arbitray $K$ corresponds to a class ( conserving the group relations
independent of $K$ ) of the MQYBE or modified quantum YB equations studied
by Gerstenhaber, Giaquinto and Schak. Various properties of the
triparametric
$\hat{R}(K;p,q)$, $\hat{R}(K;g,h)$ and $\hat{R}(K;q,h)$ are studied. In the
larger space of the modified braid equation (MBE) even $\hat{R}(K;p,q)$ can
satisfy $\hat{R}^2 = 1$ outside braid equation (BE) subspace. A
generalized, $K$- dependent, Hecke condition is satisfied by each
$3$-parameter $\hat{R}$. The role of $K$ in noncommutative geometries of
the $(K;p,q)$,$(K;g,h)$ and $(K;q,h)$ deformed planes is studied. K is
found to introduce a "soft symmetry breaking", preserving most interesting
properties and leading to new interesting ones. Further aspects to be
explored are indicated.

\end{abstract}

\vfill
\newpage

\pagestyle{plain}

\section{Introduction}
      Our starting point will be the group relations of the elements of the
quantum matrix
\begin{equation}
T = \left(\begin{array}{cc}
a & b  \\
c & d
\end{array}\right).
\end{equation}
 Three known cases will be considered:

 $(1)$ The biparametric $(p,q)$ or standard deformation of $GL(2)$

 $(2)$ The biparametric $(g,h)$ or nonstandard deformation of $Gl(2)$

 $(3)$ the $(q,h)$ or "hybrid" deformation of$GL(1/1)$.

Each set will be presented explcitly below.These three have been studied in
$[1]$ where many original sources can be found cited. We start directly
with the biparametric cases since the one - parameter deformations can then
be systematically obtained through suitable constraints ( $ p = q^{-1}$,$ g
= h$ and so on ).

    For the {\it {given group relations}} we construct for each case the
matrix $R$ satisfying
\begin{equation}
          RT_{1} T_{2} = T_{2} T_{1}R
\end{equation}

where
         $$ T_{1} = T \otimes I_{2}, T_{2} = I_{2} \otimes T$$

To start with,{\it  {we do not require}} $R$ {\it {to satisfy the Yang -
Baxter equation}} $(YB)$. It will be found that, apart from the parameters
concerned ( $ (p,q)$, $(g,h)$or $(q,h)$ ) the solution for $R$ satisfying
$(2)$ contains a supplementary arbitrary parameter $K$. Two particular
values of $K$ ( say $K_{1}$ and $K_{2}$ ) will give the two solutions of
$(YB)$ related through
\begin{equation}
         ((21)R(K_{1}))^{-1} = R(K_{2})
\end{equation}
both satisfying
                   $$ R_{12}R_{13}R_{23} - R_{23}R_{13}R_{12} = 0$$

The existence of such a pair of solutions is assured by the fact that
$(1.2)$ can be written as

        $$T_{1} T_{2}R^{-1} = R^{-1}T_{2} T_{1}$$

      The germ of our paper is in the question: what structure is obtained
when $K$ is  {\it {not}} restricted to the values $K_{1}$ and $K_{2}$.

     For arbitrary $K$ the situation is best expressed in terms of

                   $$\hat{R} = PR$$

 where $P$ is the permutation matrix and for our $4\times 4$ case it
permutes the second and the third rows.

         One obtains, for the normalizations we will choose,
\begin{equation}
            \hat{R}_{(12)}\hat{R}_{(23)}\hat{R}_{(12)} -
\hat{R}_{(23)}\hat{R}_{(12)}\hat{R}_{(23)} = \biggl( { K\over {
K_{1}}} - 1
\biggr) \biggl( {K\over { K_{2}}} - 1 \biggr) \biggl(\hat{R}_{(23)} -
\hat{R}_{(12)}\biggr)
\end{equation}
    This is our modified braid equation $(MBE)$. ( See Discussion for
further comments.)

      In terms of $R$ one obtains
\begin{equation}
  R_{12}R_{13}R_{23} - R_{23}R_{13}R_{12} = \biggl( { K\over {
K_{1}}} - 1
\biggr) \biggl( {K\over { K_{2}}} - 1 \biggr) \biggl( (123)R_{(12)} -
(213)R_{(23)}\biggr)
\end{equation}

     where $(123)$ and $(213)$ denote corresponding permutations of the
tensor factors of $V^\otimes{3} $( $R$ acting on $V\otimes V)$ . [ Having
pointed out the structure $(1.5)$ we will use throughout $(1.4)$ as our
fundamental relation.]

           Thus $(1.2)$ by itself is seen to lead to a particular class of
solutions of the "modified quantum Yang - Baxter equations" $(MQYBE)$
introduced by Gerstenhaber et al. $[2]$.Our $(1.4)$ has the same srtucture
as the equation $(2.4)$ of $[3]$ for "quantum transpositions" $(
\sigma_{12}, \sigma_{23} )$ defined by the authors , though we do {\it
{not}} impose in general their "unitarity" leading to

                        $${\hat {R}}^2 = I$$

 An example of a solution of $(1.2)$ with an arbitrary $K$ can be found in
$[1]$.

      We present below some particularly interesting explicit examples.Their
properties will reveal that the existence of such a class of more general
solutions of $(MBE)$ is more than an accident and can play a significant
role in various domains, such as noncommutative geometry.

\section{Explicit solutions}

\subsection{ Standard $(p,q)$ deformation of $GL(2)$ }

The elements $(a,b,c,d)$ of $T$ satisfy

$$ab = qba, pac = ca,$$
\begin{equation}
ad = da + (q-p)bc,\quad  pqbc = cb,
\end{equation}
$$pbd = db,\quad cd = qdc.$$

Apart from a possible normalizing factor, the solution of $(1.2)$ turns
out to be ( writing directly $\hat{R} =PR$ and assuming $p$ to be nonzero )

\begin{equation}
\hat{R}(K;p,q) =
\pmatrix{
1&0&0&0 \cr
0&(1-K)&{K\over p}&0 \cr
0&Kq&(1-{Kq\over p})&0  \cr
0&0&0&1
}.
\end{equation}

This is found to satisfy
\begin{equation}
            \hat{R}_{(12)}\hat{R}_{(23)}\hat{R}_{(12)} -
\hat{R}_{(23)}\hat{R}_{(12)}\hat{R}_{(23)} = \biggl(  K
 - 1\biggr) \biggl( {Kq\over { p}} - 1 \biggr) \biggl(\hat{R}_{(23)} -
\hat{R}_{(12)}\biggr)
\end{equation}
with
$$ K_{1} =1, K_{2} = {p\over q}$$

\subsection{ Nonstandard $(g,h)$ deformation of $GL(2)$ }

The group relations are,

$$ca =ac - gc^2, cb = bc - gdc - hac + ghc^2,$$
$$cd = dc -hc^2, da = ad - gdc + hac,$$
\begin{equation}
db = bd + g( ad - bc + hac - d^2 )
\end{equation}
$$ba =ab - h( ad - bc + hac - a^2 )$$

>From $(1.2 )$ one obtains

 \begin{equation}
\hat{R}(K;g,h) =
\pmatrix{
1&- hK&hK&ghK \cr
0&(1-K)&K&gK \cr
0&K&(1-K)& - gK  \cr
0&0&0&1
}.
\end{equation}

  This is found to satisfy
\begin{equation}
            \hat{R}_{(12)}\hat{R}_{(23)}\hat{R}_{(12)} -
\hat{R}_{(23)}\hat{R}_{(12)}\hat{R}_{(23)} = \biggl(  K
 - 1\biggr)^2 \biggl(\hat{R}_{(23)} -
\hat{R}_{(12)}\biggr)
\end{equation}

with, independently of $(g,h)$,

                   $$K_1 = K_2 = 1$$

\subsection{Hybrid $(q,h)$ deformation of $G(1/1)$}

The group relations are

$$ba = ab + hcd, ac = qca, bc = qcb,$$
\begin{equation}
dc+qcd = 0, ad - da + (1 -q)cb = 0,
\end{equation}
$$bd +db = hca,$$
$$ha^2 =hd^2 + (q +1)b^2, c^2 =0$$

One obtains from $(1.2)$

 \begin{equation}
\hat{R}(K;q,h) =
\pmatrix{
1&0&O&Kh \cr
0&(1-K)&Kq&0 \cr
0&K&(1-Kq)& 0  \cr
0&0&0&(1 - K(q+1) )
}.
\end{equation}

This satisfies
\begin{equation}
            \hat{R}_{(12)}\hat{R}_{(23)}\hat{R}_{(12)} -
\hat{R}_{(23)}\hat{R}_{(12)}\hat{R}_{(23)} = \biggl(  K
 - 1\biggr)\biggl( Kq - 1 \biggr) \biggl(\hat{R}_{(23)} -
\hat{R}_{(12)}\biggr)
\end{equation}
Here $K_1 = 1, K_2 = q^{-1}$, both being independent of $h$.

\section{ Properties}

\subsection { K and triangularity }

 The matrix $R$ is called "triangular" if
$$(21)R =R^{-1}$$  when $$\hat{ R}^2 = (PR)^2 = I$$

In $[2]$ the term "unitary" is used in this context.
For an $R$-matrix satisfying the Yang -Baxter constraint $(YB)$ the
following features are well known.

 (1): For standard($q$ or$(p,q)$ ) deformations the $R$ satisfying $(YB)$ is
"quasitriangular" and

               $$ \hat{ R}^2  \neq I$$

(2): For nonstandard ( $h$ or $(g,h)$ ) deformations for $R$ satisfying
$(YB)$ one has "triangularity" or
                     $$\hat{ R}^2 = I$$

( It is in this sense that we use the term triangular, without $R$ being
necessarily strictly upper or lower triangular.)

    In presence of an arbitrary $K$ the modified braid equation $(MBE)$
{\it {breaks this dichotomy}}. Specificaly in the preceding three cases one
has the following situation:
\begin{equation}
         (21) R(K;p,q)  = ( R(K';p,q) )^{-1}
\end{equation}

where
             $$ K' = K ( K(1+q{ p}^{-1}) - 1)^{-1}$$

\begin{equation}
         (21) R(K;g,h)  = ( R(K';g,h) )^{-1}
\end{equation}

where
             $$ K' = K (2K - 1)^{-1}$$

 \begin{equation}
         (21) R(K;q,h)  = ( R(K';q,h) )^{-1}
\end{equation}

where
             $$ K' = K ( K(1+q) - 1)^{-1}$$

 Thus in each case one obtains

\begin {equation}
K' = K\biggl({K\over { K_{1}}} + {K\over { K_{2}}} - 1 \biggr)^{-1}
\end{equation}

In general none is triangular (or unitary). On the other hand { \it {in
each case one can have triangularity }} by choosing
$$K'=K$$ or
\begin{equation}
K = 2K_{1}K_{2}(K_{1}+K_{2})^{-1}
\end{equation}

For the three previous cases this gives respectively
\begin{equation}
K = 2p(p+q)^{-1},\qquad 1, \qquad 2(1+q)^{-1}
\end{equation}

Thus for the nonstandard case triagularity coincides with the $(YB)$
property. In contrast, for the other two cases triangularity implies a
nonzero right hand side in $(1.4)$. In particular for the $(p,q)$ case one
obtains (permuting the second and the third rows of $\hat {R}$) for

$$K = 2p(p+q)^{-1}$$

\begin{equation}
R =\pmatrix{
1&0&0&0 \cr
0&{2pq}\over {p+q}&{p-q}\over {p+q}&0 \cr
0&{q-p}\over {p+q}&2\over {p+q}&0  \cr
0&0&0&1
}.
\end{equation}

Now one has,

$$ R(K;p,q) = ((21)M)^{-1}M$$

where one can set,choosing an upper triangular form,

\begin{equation}
M=\pmatrix{
1&0&0&0 \cr
0&({{2pq}\over {p+q}})^{1\over 2}&{p-q}\over {(2pq(p+q))^{1\over 2}}&0 \cr
0&0&({{p+q}\over {2pq}})^{1\over 2}&0  \cr
0&0&0&1
}.
\end{equation}

 $R$ is invariant under $$M\rightarrow {VM}$$ where $$(21)V = V$$.

\subsection {Projectors}

For each case ( $ (K;p,q),(K;g,h),(K;q,h)$) one obtains , $I$ being the
$4\times 4$ unit matrix,

\begin{equation}
\hat{R}^2 = X \hat{R} + (1 - X)I
\end{equation}
$$ X = 2 -\biggl({K\over { K_{1}}} + {K\over { K_{2}}} \biggr) $$

 Thus for the three cases $((p,q),(g,h),q,h))$ one has respectively

$$X= 2 - K(1 + qp^{-1})$$
\begin{equation}
X = 2(1 - K)
\end{equation}
$$X = 2 - K(1+q)$$

Two special cases are

$$ X=0, \qquad (\hat{R})^2 = I$$
and
$$ X = 2, \qquad (\hat {R} - I)^2 = 0$$

For $X\neq 2$ one obtains for each deformation considered two projectors
$(P^2 =P)$ as follows

\begin{equation}
P_{1} = {{(\hat{R} - I )}\over {(X -2)}}
\end{equation}
\begin{equation}
P_{2} = {{(\hat{R} -(X-1)I)}
\over {(2 -X)}}
\end{equation}

Finally one has

\begin{equation}
\hat R = (X-1)P_{1}+P_{2}
\end{equation}
with
$$P_{1}+P_{2}= I, \qquad P_{1}P_{2}=0$$

Note that one obtains a canonical formalism valid for all the deformations
considered.

It follows from the preceding results that if $C$ is a column vector (
with $4$ rows) a constraint ( $n$ being a constant)

                      $$C=n\hat {R}C$$
is only consistent with the actions of the projectors for

$$n=1; \qquad (P_{1}C=0, C=P_{2}C)$$
or

$$n=(X-1)^{-1}; \qquad (P_{2}C=0, C=P_{1}C)$$

This fact should be kept in mind for what follows.

\section{K and noncommutative planes}

Detailed study of noncommutativity implemented via $\hat{R}$ in the plane
and higher dimensional spaces can be found in $[4,5,6]$ where numerous
sources are cited. Here we limit our considerations to the two dimensional
plane. But we let our $\hat {R}$ be more general by letting it depend on an
exrta arbitray parameter $K$ permitted by our $(MBE)$. Our $\hat {R}$ will
depend on $3$ parametres. The biparametric nonstandard deformation with
differential calculus was first presented ( for $K=1, g=h, h=h'$ )
in $[7]$. The original formalism is due to Wess and Zumino $[8]$.

  We use the following notations:

    $$ x^{i}=(x^{1},x^{2}) = (x,y)$$
$$ dx^{i}={\xi}^{i}=({\xi}^{1},{\xi}^{2})= (\xi,\eta)$$
$$ (\alpha,\beta) = (p,q),(g,h),(q,h)$$

 We postulate
\begin{equation}
x^{i}x^{j}=\bigl (\hat{R}(K;\alpha,\beta)
\bigr)^{ij}_{i'j'}x^{i'} x^{j'}
\end{equation}
i,e,
$$ \bigl ( P_{1} \bigr )^{(ij)}_{i'j'}x^{i'}
x^{j'} = 0$$

\begin{equation}
{\xi}^{i}{\xi}^{j}= - {1\over {(1-X)}} \bigl (\hat{R}(K;\alpha,\beta)
\bigr)^{(ij)}_{i',j'}{\xi}^{i'} {\xi}^{j'}
\end{equation}
i,e,
$$ \bigl ( P_{2} \bigr )^{(ij)}_{{i'},{j'}}{\xi}^{i'}
{\xi}^{j'} = 0$$

\begin{equation}
x^{i}{\xi}^{j}={1\over {(1-X)}}\bigl (\hat{R}(K;\alpha,\beta)
\bigr)^{(ij)}_{{i'},{j'}}{\xi}^{i'} x^{j'}
\end{equation}
where
$$(1-X)=\biggl({K\over { K_{1}}} + {K\over { K_{2}}} -1 \biggr) $$

The bilinear constraints (4.28),(4.29),(4.30), related through
derivations, are required to satisfy  suitable consistency relations.
(See,for example, Sec.$4$ of $[4]$ and $[7]$.) Following the usual procedure
the required consistency for our case can be shown to be assured precisely
by our generalized Hecke condition, namely,

\begin{equation}
\biggl (\hat {R}(K;\alpha,\beta) -I\biggr )\biggl ({\hat
{R}({K;\alpha,\beta)}\over {(1-X)}} +I\biggr ) =0
\end{equation}

or $$P_{1}P_{2}=0$$

This generalizes some well known results. Thus,for example, setting
     $$p=q^{-1}, \quad K_{1}=1 \quad K= K_{2} = q^{-2}$$
and changing the normalization of $R$ by a factor $q$ one obtains the
result $(4.4.15)$ of $[4]$. One obtains analogous generalizations for the
other cases. Note that the consistency is obtained for our case by {\it
{by implementing }}$K$ {\it {nontrivially}} through the factor $(1-X)$ for
the $(\xi,\eta)$ constraints. But once this is done the final
consequences of $(4.28)$ and $(4.29)$ turn out to be {\it {systematically
independent of}} K. ( Those of $(4.30)$ do involve $K$ but, as will be
shown below, in a particularly simple fashion.) We recapitulate for
completeness the first two sets of results which are the same as one
would obtain with $ K = (K_{1},K_{2})$.

One obtains for $(\alpha,\beta) = (p,q)$

\begin{equation}
pxy = yx
\end{equation}
\begin{equation}
{\xi}^{2}=0, \qquad {\eta}^{2}=0, \qquad  \eta \xi +q\xi \eta = 0
\end{equation}

For $(\alpha,\beta) = (g,h)$

\begin{equation}
xy - yx = gy^{2}
\end{equation}
\begin{equation}
{\xi}^{2}= h\xi \eta , \qquad {\eta}^{2}=0, \qquad  \eta \xi +\xi \eta = 0
\end{equation}

The results above are for $GL(2)$. For $GL(1/1)$, namely for
$$(\alpha,\beta) = (q,h)$$
one obtains

\begin{equation}
xy = qyx, \qquad y^{2} =0
\end{equation}
\begin{equation}
(1 + q){\xi}^{2} +h{\eta}^{2}=0,  \qquad  \eta \xi +\xi \eta = 0
\end{equation}

For the deformed $GL(1/1)$ $y$ becomes fermionic. ( After exhibiting as
above how the three cases can be treated uniformly in our formalism, in what
follows we will consider only the deformations $(K;p,q)$ and $(K;g,h)$ of
$GL(2)$. Those for $GL(1/1)$ can easily be added.)

In contrast to the foregoing results, the consequences of $(4.30)$ {\it
{involve}} K {\it {nontrivially}}. For $\hat {R}(K;p,q)$ one obtains

$$ x{\xi} = {1\over {(1-X)}}{\xi}x, \qquad x{\eta} = {1\over
{(1-X)}}\biggl({\xi}y + {K\over {p}} {\Phi}_{1}\biggr)$$
\begin{equation}
y{\xi} = {1\over {(1-X)}}\biggl({\eta}x - { Kq\over {p}} {\Phi}_{1}\biggr),
\qquad y{\eta} = {1\over {(1-X)}}{\eta}y
\end{equation}
with
\begin{equation}
         {\Phi}_{1} = ({\eta}x - p{\xi}y )
\end{equation}
For $p=q^{-1}$, $K= q^{-2}$ and again suitably choosing the normalization
of $\hat R$ these results reduce to $(4.1.8)$ of $[4]$. In order to
compare with $\kappa$ of $(4.1.10)$ of $[4]$ one can show by reordering
terms
 \begin{equation}
         {\Phi}_{1}^{2} = {1\over {(1-X)}}(-qp + Kqp +
p^{2}- Kqp + pq - p^{2})({\xi}{\eta}xy) = 0
\end{equation}

Note that $K$ reappears on reordering but the coefficient of $K$ in
the numerator vanishes separately. Thus, apart from the overall factor, $K$
{\it {appears as a factor of the nilpotent}}
${\Phi}_{1}$. Moreover one can show that

$$ px{\Phi}_{1} ={1\over {(1-X)}} K{\Phi}_{1}x,\qquad   y{\Phi}_{1} =
{1\over {(1-X)}}Kq{\Phi}_{1}y$$
\begin{equation}
{1\over {(1-X)}}(p+q- Kq){\xi}{\Phi}_{1} = -{\Phi}_{1}\xi,\qquad  {1\over
{(1-X)}}(p+q- Kq)\eta{\Phi}_{1} = -pq{\Phi}_{1}\eta
\end{equation}

For the prescriptions indicated before ($K=q^{-2}$ and so on) one finds
back the corresponding results of Sec.$4.1.13$ and Sec.$4.1.17$ of [4].

For $\hat {R}(K;g,h)$ one obtains (compare $(3.1.6)$ of $[5]$ where
$K=1$ and $g=h$)

$$ x{\xi} = {1\over {(1-X)}}\biggl({\xi}x + Kh\Phi_{2}\biggr), \qquad
x{\eta} = {1\over {(1-X)}}\biggl({\xi}y + K{\Phi}_{2}\biggr)$$
\begin{equation}
y{\xi} = {1\over {(1-X)}}\biggl({\eta}x -  K {\Phi}_{2}\biggr), \qquad
y{\eta} = {1\over {(1-X)}}{\eta}y
\end{equation}
with
\begin{equation}
         {\Phi}_{2} = ({\eta}x - {\xi}y + g{\eta}y)
\end{equation}
and
\begin{equation}
         {\Phi}_{2}^{2} = 0
\end{equation}

Moreover,

$$ x{\Phi}_{2} ={1\over {(1-X)}}\biggl( K{\Phi}_{2}x + K(g
-h){\Phi}_{2}y\biggr) ,\qquad   y{\Phi}_{2} ={1\over {(1-X)}}
Kq{\Phi}_{2}y$$
\begin{equation}
{1\over {(1-X)}}(2- K){\xi}{\Phi}_{2} = -({\Phi}_{2}\xi+
(h-g){\Phi}_{2}\eta),\qquad       {1\over {(1-X)}}(2- K)\eta{\Phi}_{2} =
-{\Phi}_{2}\eta
\end{equation}

The results for the$(g,h)$ case can of course be obtained independently. But
they are obtained more efficiently and with a deeper understanding by
starting from the corresponding ones for $(p,q)$ and using the
"contraction" studied in the following section. It is instructive to see,
in particular, how the $(g-h)$ factors in the results above arise (end of
the next section). These terms are present even for the $(YB)$ subspace
($K=1$) unless $g=h$. Finally, for $K=1$ and $g=h$ one obtains the simple
results of Sec.$4.1.17$ of $[4]$.

\section { $((K;p,q) \rightarrow (K;g,h))$ : singular limit of a
transformation }

In Sec.$4$ of $[1]$ such a passage was presented for the case where
$R(p,q)$ and $R(g,h)$ both satisfied ($YBE$). Here we generalize it to
include an arbitrary $K$. In fact the same transformation will work
again, leading to a well defined $R(K;g,h)$. We want to emphasize this
fact. It underlines again the "soft symmetry breaking" role of $K$.
Moreover we will display here how the corresponding features of the two
noncommutative plains are related systematically through this
"contraction"procedure. The $(K;g,h)$ -deformed plane emerges in full
detail from the $(K;p,q)$-deformed one.  Some previous sources are cited in
$[1]$, which in turn lead to some original ones.

Setting
\begin{equation}
G = \left(\begin{array}{cc}
1 & \omega  \\
0 & 1
\end{array}\right).
\end{equation}
and with $R = P\hat {R}$ one obtains

 \begin{equation}
(G^{-1}\otimes G^{-1})R(K;p,q)(G\otimes G) =
\pmatrix{
1&- K(q-1)\omega & {K\over {p}}(q-1)\omega & - {K\over {p}}
(p-1)(q-1){\omega}^2
\cr 0 & Kq & (1-K{q\over {p}}) & K{q\over {p}} (p-1){\omega} \cr 0
& (1-K) & K\over {p} &
 {K\over {p}} (p-1){\omega}\cr 0&0&0&1
}.
\end{equation}

Now, as in $[1]$, let $p\rightarrow 1, q\rightarrow 1$ in such a
way that $(p-1){(q-1)}^{-1}$ remains constant. And ${\omega}_{0}$ being a
constant, define
$$\omega = {\omega}_{0}((p-1)(q-1))^{-{1\over 2}}$$

Now one can define finite constants $(g,h)$ such that as $p\rightarrow 1$
and $q\rightarrow 1$

\begin{equation}
 ((1-p)\omega )\rightarrow g, \qquad ((q-1)\omega ) \rightarrow h
\end{equation}

Now from $(2.10)$ and $(5.47)$ ( with $R = P\hat {R} )$ ), one obtains
\begin{equation}
 (G^{-1}\otimes G^{-1})R(K;p,q)(G\otimes G) \rightarrow R(K;g,h)
\end{equation}

The same procedure works for $(a,b,c,d)$, $(x,y)$ and $(\xi,\eta)$. In this
section, to distinguish the cases $(p,q)$ and $(g,h)$, we will use for the
latter the notations
$$(\tilde{a},\tilde{b},\tilde{c},\tilde{d});(\tilde{x},\tilde{y});
(\tilde{\xi},\tilde{\eta})$$

Consistently with transformation of $R$ one defines ( with the previous
defini!tions of $G$ and $T$ )

$$G^{-1}TG = \left(\begin{array}{cc}
           \tilde{a} & \tilde{b}  \\
           \tilde{c} & \tilde{d}
            \end{array}\right).$$

$$G^{-1}\left(\begin{array}{c}
           x  \\
           y
        \end{array}\right) =\left(\begin{array}{c}
                                \tilde{x}  \\
                                \tilde{y}
                              \end{array}\right).$$

$$G^{-1}\left(\begin{array}{c}
          \xi \\
          \eta
        \end{array}\right) =\left(\begin{array}{c}
                                \tilde{\xi}  \\
                                \tilde{\eta}
                              \end{array}\right).$$

Let us now consider some examples to appreciate how the technique works.
>From the preceding definitions one obtains

$$ \tilde{a}=a-{\omega}c, \qquad
\tilde{b}=(b-{\omega}d)+{\omega}(a-{\omega}c)$$
\begin{equation}
\tilde{c}=c, \qquad \tilde{d}=d+{\omega}c
\end{equation}

The inverse relations are easily obtained. Using  them and the group
relations for $(a,b,c,d)$ one obtains
$$
\tilde{c}\tilde{a} = c(a - {\omega}c) = pac - {\omega}c^2 =
p(\tilde{a}+{\omega}\tilde{c})\tilde{c} - {\omega}{\tilde{c}}^{2} =
p\tilde{a}\tilde{c} - (1 -p) {\omega}{\tilde{c}}^{2}
$$

Using the definition of $g$ now one obtains, in the limit,
\begin{equation}
 \tilde{c}\tilde{a}=\tilde{a}\tilde{c} - g{\tilde{c}}^2
\end{equation}

Again,
$$
\tilde{c}\tilde{b}=c({\omega}a + b - {\omega}^{2}c - {\omega}d )
                  =pqbc - {\omega}qdc + p{\omega}ac - {\omega}^{2}c^2$$
                $$  = pq\tilde{b}\tilde{c} +q(p-1)
{\omega}\tilde{d}\tilde{c} - p(q-1){\omega}\tilde{a}\tilde{c} +(1-p)(q-1)
{\omega}^{2}{\tilde{c}}^2
$$
giving in the limit
\begin{equation}
 \tilde{c}\tilde{b}
                  = \tilde{b}\tilde{c} -g\tilde{d}\tilde{c}
-h\tilde{a}\tilde{c} +gh{\tilde{c}}^2
\end{equation}

We have thus obtained the first two group relations ( with tildes added to
avoid confusion) for the nonstandard case $(2.9)$. The others can be
obtained analogously. Let us now look at the $(K;g,h)$-deformed plane. One
obtains from the definitions introduced
\begin{equation}
\tilde{x}= x - {\omega}y, \qquad \tilde{y} =y
\end{equation}

Hence, using the constraints for $(x,y)$
$$
\tilde{x}\tilde{y}-\tilde{y}\tilde{x}=(x -{\omega}y)y - y(x -{\omega}y)
=xy -yx= (1-p)xy = (1-p)(\tilde{x}\tilde{y}+{\omega}y^2)
$$
Now taking limit and using the definition of $g$,
\begin{equation}
\tilde{x}\tilde{y}-\tilde{y}\tilde{x}= g {\tilde{y}}^2
\end{equation}
This is the nonstandard version (with tildes added)

Similarly starting with
\begin{equation}
\tilde{\xi}  = \xi - \omega \eta, \tilde{\eta} = \eta
\end{equation}
and using
$$ {\xi}^2 =0,\qquad {\eta}^2 =0,\qquad ({\xi}{\eta}+q{\eta}{\xi})=0$$
one obtains in the limit the expected results
$$
 {\tilde{\xi}}^2 = h\tilde{\xi}\tilde{\eta}, \qquad
{\tilde{\eta}}^2=0, \qquad
(\tilde{\xi}\tilde{\eta}+\tilde{\eta}\tilde{\xi})=0
$$

These simple cases have been presented to give a feeling for the limiting
process at work. But they have further usefulness. For the important
nilpotent operators of the preceding section one easily obtains, taking our
limits,

 $${\Phi}_{1} =({\eta}x -p{\xi}y) \rightarrow
(\tilde{\eta}\tilde{x} - \tilde{\xi}\tilde{y} + g \tilde{\eta}\tilde{y})
\quad  = {\Phi}_{2}$$

Hence avoiding a lengthy reordering process one obtains directly from

$${{\Phi}_{1}}^2 = 0$$

$$ ({{\Phi}_{1}}^2) \rightarrow {{\Phi}_{2}}^2 = 0$$

The commutators of ${\Phi}_{2}$ can again be derived simply from those of
 ${\Phi}_{1}$. The terms involving $(g -h)$ in the latter set can be seen
to emerge as follows:
$$
  (p^{-1} - q)\omega  =  ( p^{-1} (1-p)\omega - (q-1)\omega)
\rightarrow (g -h)
$$

$$
  (pq - 1)\omega  =  ( p(q-1)\omega - (1-p)\omega)
\rightarrow (h -g)
$$

 \section{Discussion}

The following points are worth noting.

$(1)$: If $\hat{R}(K;\alpha,\beta)$ depend linearly on $K$ and satisfy the
braid equation for $K=(K_{1},K_{2})$, then the right hand side of $(1.4)$
becomes almost evident as follows. One can set
\begin{equation}
            \hat{R}_{(12)}\hat{R}_{(23)}\hat{R}_{(12)} -
\hat{R}_{(23)}\hat{R}_{(12)}\hat{R}_{(23)} = \biggl( { K\over {
K_{1}}} - 1\biggr) \biggl( {K\over { K_{2}}} - 1 \biggr) Z
\end{equation}

The first two factors on the right assure the braid property for
$K=(K_{1},K_{2})$. Next one notes the follwing points:

The left hand side is trilinear in $K$. Hence $Z$, coming after the first
two factors, should be linear in $K$.

The left side is antisymmetric under the exchange
$$ \hat{R}_{(12)} \leftrightarrow   \hat{R}_{(23)}$$
Hence $Z$ should have the same property.Thus the evident ansatz is
\begin{equation}
Z = \biggl(\hat{R}_{(23)} -\hat{R}_{(12)}\biggr)
\end{equation}

This is indeed found to be correct. A possible $K$-independent constant
factor can be normalized to unity, as we have done.

 (The following two properties have been pointed out to the author by Daniel
Arnaudon.)

  $(2)$: For the $(p,q)$ and the $(q,h)$ cases one can write
\begin{equation}
 \hat{R}(K;\alpha,\beta) = c_{1} \hat{R}(K_{1};\alpha,\beta)
+c_{2}\hat{R}(K_{2};\alpha,\beta)
\end{equation}
where
$$ c_{1}+c_{2}=1;\qquad   c_{1}K_{1}+c_{2}K_{2}=K$$

However, for the $(g,h)$ case (since $ K_{1}=K_{2}=1$) such a relation does
not hold for $ K \neq 1$.

$(3)$: For
\begin{equation}
\hat{R}_{(12)}\hat{R}_{(23)}\hat{R}_{(12)} -
\hat{R}_{(23)}\hat{R}_{(12)}\hat{R}_{(23)} = \lambda \biggl(\hat{R}_{(23)} -
\hat{R}_{(12)}\biggr)
\end{equation}
and
\begin{equation}
\hat{R}^2 = X \hat{R} + (1 - X)I
\end{equation}

( $\lambda$ and $X$ not being necessarily restricted to the values
considered previously) defining
\begin{equation}
\hat{S} = ( \hat{R} - {\mu}I )
\end{equation}
one can verify that

\begin{equation}
\hat{S}_{(12)}\hat{S}_{(23)}\hat{S}_{(12)} -
\hat{S}_{(23)}\hat{S}_{(12)}\hat{S}_{(23)} = (\lambda +X\mu -{\mu}^2)
\biggl(\hat{S}_{(23)} -\hat{S}_{(12)}\biggr)
\end{equation}
This generalizes an analogous result of $[2]$ since we do not restrict $R$
to be "unitary". One can choose $\mu$ so that $\hat{S}$ satifies the
braid equation.Directly connected with the last two equations is the
following, canonical relation valid for all the cases considered before,

\begin{equation}
 \hat{R}(K;\alpha,\beta) = \biggl( {K\over { K_{i}}}
\biggr)\hat{R}(K_{i};\alpha,\beta) -\biggl( {K\over { K_{i}}} -1
\biggr)I
\end{equation}
Here $K_{i}$ denotes either one of the "braid" ( or YB) values of K. This
is of basic importance. The parameter $K/(K_{i}-K)$ can be shown to provide
the prescription for Baxterization. In fact, $(6.63)$ can be recognized to
correspond to the usual ansatz for Baxtarization [9].

$(4)$: The works of Gerstenhaber, Giaquinto and Schak $[2,3]$ assure
that our $(MBE)$ encodes deformations satisfying basic criteia
but removing certain restrictive features of the standard $(BE)$ (
or$(YB)$). For our class the factor $\lambda$ on the right in $(6.59)$ is
neither zero nor entirely arbitrary. It has the specific form given by
$(1.4)$ arising out of our basic condotion: $K$-independence of the group
relations.Our parametrization of of this factor carries information. The YB
or the braid solutions are obtained effortlessly as byproducts. This leads
also to agreeable properties dsignated here as "soft symmetry breaking "
role of
$K$ in the noncommutative geometries studied. For all $K$ (and all the case
studied) one obtains the crucial, canonical Hecke condition we have
emphasized. It permits us to introduce consistently and uniformly the
noncommutativity constraints. Let us recapitulate the remarkable
consequences.

$(a)$: The bilinear constraits for the coordinates  and those for the
differntials remain {\it {independent}} of $K$.

$(b)$: In the constraints involving both coordinates and differentials $K$
does appear but in a "minimal" fashion. Apart from a simple overall factor,
in the linear combinations on the right $K$ appears as a factor of a
nilpotent combination$({\Phi}_{1}$ or
${\Phi}_{2}$). Along with the commutation relations satisfied by these
nilpotents, this has the consequence that reordering any higher order
product one obtains, apart from an overall factor, finally {\it {linearity}}
in $K$.The operator $\Phi$, crucial for constructing covariant derivatives
$[4]$, remains nilpotent for arbitrary $K$ (see $(4.40)$). The main point
is that conserving the $(x,y)$ and the $(\xi,\eta)$ commutators and without
violating the constraints imposed by the postulated actions of exterior
derivations one can implement the parameter $K$ in the mixed commutators
($(x,\xi)$,etc.), even there conserving good properties.

$(c)$: The "contraction" procedure leading from standard ($p,q$) to
nonstandard ($g,h$) deformations is not perturbed by $K$. Even the titles
of previous papers $[10,11]$ give an idea of the scope of this approach. It
is resuring to note that one can continue to implement it in presence of
$K$.

Having noted some interesting features of the results obtained let us now
look at further developments they suggest. One naturally thinks of the
following aspects:

$(1)$: Extension of our results to higher dimensional $(MBE)$. Firstly by
going beyond the $4\times 4$ cases for deformations of $GL(2)$. Secondly
by starting from group relations for deformed $SL(N)$ and $SO(N)$. Higher
dimensional cases have already been studied in $[2,3]$. Our aim would be
to obtain explicit srtuctures correponding to {\it {conserved group
relations}} for such cases. Then one can see if our soft symmetry breaking
still gets implemented and in what fashion.

$(2)$: For $K= (K_{1},K_{2})$ the $R$ matrix flips the tensor components of
coproducts. Having obtained more general modified $R$ mtrices it would be
important to examine the consequences for coproducts as $K$ moves away
from the $(YB)$ values.

$(3)$: A more complete study of the role of $K$ in noncommutative geometries
induced by $3$-parameter deformations $(K;\alpha,\beta)$. Even for the
$2$-plane we have stopped at a certain point leaving much to be done.
After constructing higher dimensional matrices $\hat{R}(K;\alpha,\beta)$
one can implement them in higher dimensional spaces.

$(4)$: Study of twists in the context of "modified" $R$ matrices. In
particular, the fact that one can implement triangularity for all types of
deformations by suitablly choosing $K$ suggests intriguing possibilities.
Various aspects studied in $[12,13,14]$ can be reexamined in this broader
context.

$(5)$: Our MBE ( or MYBQE of Gerstenhaber et al.) and Baxterization can be
seen to be ( see $(6.63)$ ) two facets of the same underlying construction,
namely the general solution of the $RTT$ relations. In the first case the
parameter $K$ is kept fixed in each term and the right hand side of the
braid equation is allowed to be nonzero. In the second one the right is
held fixed at zero and to permit this the parameter is suitably varied
from term to term. The two procedures are complementary! This links $MBE$
with integrable models.

$(6)$:  What are the consequences for knot invariants associated to an
$\hat {R}$ as $K$ moves away from the "braid values"? Can a conceptually
consistent generalization ( parametrized third Reidemeister move) be
implemented fruitfully ?

Presumably this list is not exhaustive. Some of these objectives should be
directly accessible. Elsewhere one may encounter obstructions. We hope to
explore different directions in future studies.

\smallskip

\smallskip

This work owes much to sustained and reassuring help from Daniel Arnaudon.
It goes beyond results explicitly attributed to him. Our treatment of
noncommutative planes took shape from succesive discussions with John
Madore.

\smallskip

\smallskip

\bibliographystyle{amsplain}

\end{document}